\documentclass[a4paper]{article}
\usepackage{amsmath,amssymb}
\usepackage[latin1]{inputenc}
\usepackage[dvips]{graphics}
\input{diagrams.tex}

\newcommand{\NN}{\mathbb{N}}
\newcommand{\ZZ}{\mathbb{Z}}
\newcommand{\CC}{\mathbb{C}}

\newcommand{\Mod}{\operatorname{Mod}}
\newcommand{\Ind}{\operatorname{Ind}}
\newcommand{\Supp}{\operatorname{Supp}}
\newcommand{\Hom}{\operatorname{Hom}}
\newcommand{\Ext}{\operatorname{Ext}}

\newcommand{\QT}{{Q}_T}
\newcommand{\QC}{{Q}_C}

\newcommand{\CatT}{\mathcal{C}_T}
\newcommand{\CatC}{\mathcal{C}_C}

\newcommand{\rel}{\operatorname{Rel}}
\newcommand{\eqv}{\Theta}

\newcommand{\idc}{\operatorname{id}_\CC}

\newcommand{\ot}{\leftarrow}
\newcommand{\za}{\alpha}
\newcommand{\zb}{\beta}

\newcommand{\zg}{\gamma}

\newtheorem{theorem}{Theorem}[section]
\newtheorem{proposition}[theorem]{Proposition}
\newtheorem{conjecture}[theorem]{Conjecture}
\newtheorem{corollary}[theorem]{Corollary}
\newtheorem{lemma}[theorem]{Lemma}

\newtheorem{definition}{Definition}

\newtheorem{remark}[theorem]{Remark}

\newenvironment{proof}{\begin{trivlist}\item{\bf{Proof.}}}
  {\hfill\rule{2mm}{2mm}\end{trivlist}}

\title{Quivers with relations arising from clusters \\ {(}$A_n$ case{)}}
\author{P. Caldero, F. Chapoton, R. Schiffler}
\date{\today}
\begin{document}
\maketitle
\begin{abstract}
  Cluster algebras were introduced by S. Fomin and A. Zelevinsky in
  connection with dual canonical bases. Let $U$ be a cluster algebra
  of type $A_n$. We associate to each cluster $C$ of $U$ an abelian
  category $\CatC$ such that the indecomposable objects of $\CatC$ are
  in natural correspondence with the cluster variables of $U$ which
  are not in $C$. We give an algebraic realization and a geometric
  realization of $\CatC$. Then, we generalize the ``denominator
  Theorem'' of Fomin and Zelevinsky to any cluster.
\end{abstract}

\setcounter{section}{-1}

\section{Introduction}

Cluster algebras were introduced in the work of A. Berenstein, S.
Fomin, and A. Zelevinsky, \cite{cluster1,cluster2,cluster3,Ysystems}.
This theory appeared in the context of dual canonical basis and more
particularly in the study of the Berenstein-Zelevinsky conjecture.
Cluster algebras are now connected with many topics: double Bruhat
cells, Poisson varieties, total positivity, Teichmüller spaces. The
main results on cluster algebras are on the one hand the
classification of finite cluster algebras by root systems and on the
other hand the realization of algebras of regular functions on double
Bruhat cells in terms of cluster algebras.

Recall some facts about cluster algebras. Cluster algebras $U$ of rank
$n$ form a class of algebras defined axiomatically in terms of a
distinguished set of generators $\{u_1,\ldots, u_n\}$. A cluster is a
set of ``cluster variables'' $\{w_1,\ldots, w_n\}$ obtained
combinatorially from $\{u_1,\ldots, u_n\}$. The so-called Laurent
phenomenon asserts that each cluster variable is a Laurent polynomial
in the set of variables given by a cluster. For each cluster $C$, one
can define combinatorially an oriented quiver $\QC$.

Suppose from now on that $U$ is a finite cluster algebra,
\textit{i.e.}  there only exists a finite number of cluster variables.
It is known that $U$ can be described by the data of a root system
$X_n$. Moreover, there exists a cluster $\Sigma$ such that $Q_\Sigma$
is the alternated quiver on $X_n$. S. Fomin and A. Zelevinsky give in
\cite{cluster2} a more precise description of the Laurent phenomenon:
there exists a one-to-one correspondence $\alpha\mapsto w_\alpha$
between the set of almost positive roots of $X_n$, \textit{i.e.}
positive roots and simple negative roots, and the set of cluster
variables, such that the denominator of $w_\alpha$ as a Laurent
polynomial in $\Sigma$ is given by the decomposition of $\alpha$ in
the basis of simple roots. Via the Gabriel Theorem, this property
suggests a link between cluster algebras and representation theory of
artinian rings. This is what investigate the authors of \cite{MRZ},
\cite{marall}.

In this paper, we give conjectural relations $R_C$ on the quiver
$\QC$, such that for any cluster $C$, the denominators of the cluster
variables as Laurent polynomial in $C$ are described by
indecomposables of the category $\CatC$ of representations of $\QC$
with relations $R_C$.

The main result of this article is the proof of this conjecture in the
$A_n$ case.

Another important result of this paper is a geometric realization of
the ca\-tegory $\CatC$ in the $A_n$ case. Recall that the algebra of
regular functions on the 2-grassmannian of $\CC^{n+3}$ is a finite
cluster algebra $U$ of type $A_n$. Via this realization, the cluster
variables are in natural bijection with the diagonals of a regular
$(n+3)$ polygon. Moreover, a result of Fomin and Zelevinsky asserts
that this bijection gives a one-to-one correspondence between the set
of clusters of $U$ and the set of diagonal triangulations of the
polygon. Theorem \ref{equivalence} gives a simple realization of the
category $\CatC$ in terms of the diagonals of the $(n+3)$ polygon.
There also exists a more canonical category associated to a finite
cluster algebra and also studied in \cite{marall}. We give in the
$A_n$ case a geometric realization of this category, see Theorem
\ref{orbit}.

\section{Quivers of cluster type}

\label{sectionQCT}

Let $X_n$ be a simply-laced Dynkin diagram of rank $n$ and finite
type. First we need to recall some material on clusters.

Each cluster $C$ of a cluster algebra of rank $n$ is associated with a
sign-skew-symmetric square matrix $B_C$ whose lines and columns are
indexed by the cluster variables of the cluster $C$. In the
simply-laced case, coefficients of the matrices $B_C$ belong to the
set $\{-1,0,1\}$. Hence it is easy and convenient to depict these
matrices using oriented graphs (once a convention is chosen for the
orientation). This oriented graph is called the quiver associated to
the cluster $C$ and is denoted by $\QC$. It is known that all
triangles (and more generally cycles) in these quivers are oriented in
a cyclic way \cite[Proposition 9.7]{cluster2}.

The mutation procedure of clusters contains in particular a mutation
rule for the associated matrices, which can be translated as a
mutation rule for the associated quivers. In the simply-laced case,
the mutation rule can be further simplified. The result is as follows.

Let $C$ be a cluster in a cluster algebra of simply-laced type. The
mutation $\mu_i(\QC)$ of the quiver $\QC$ at a vertex $i$ is described
as follows. First, all arrows incident to $i$ in $\QC$ are reversed in
the mutated quiver. Then, for each pair of one incoming arrow $j\to i$
and one outgoing arrow $i \to k$ in $\QC$, the arrow $j \to k$ is in
the mutated quiver if and only if the arrow $k \to j$ is not in $\QC$.
The other arrows of $\QC$ are kept unchanged in the mutated quiver.

By definition, a \textit{shortest path} in the quiver $Q_C$ is an
oriented path (with no repeated arrow) contained in an induced
subgraph of $Q_C$ which is a cycle.

\begin{definition}
  \label{conject-relation}
  Let $\QC$ be the quiver associated to the cluster $C$. For each
  arrow $ i \to j$ in $\QC$, a relation $\rel_{i,j}$ is defined as
  follows. Consider the set of shortest paths from $j$ to $i$:
  \begin{itemize}
  \item If there are exactly two distinct paths $c$ and $c'$ then $c=c'$.
  \item If there exists only one path $c$ then $c=0$.
  \item If there is no such path, there is no relation.
  \end{itemize}

  To each cluster $C$, one defines the (abelian) category $\Mod \QC$ of
  representations of the quiver $\QC$ modulo the relations $\rel_{i,j}$
  for all arrows $i\rightarrow j$ of $\QC$.
\end{definition}

Remark: it is clear in type $A$ that there is at most one shortest
path for each arrow. One could show using the geometric model of type
$D$ cluster algebras that there are at most two such paths in this
case. In order for the conjecture to make sense also in the
exceptional cases $E$, it remains to prove that this is also true in
these cases. We will not consider this question here.

Obvious remark: there is a natural one-to-one correspondence between
the vertices of the quiver and the simples of $\Mod \QC$. Therefore, the
isomorphism class of the simple module associated to the vertex $i$
will be denoted $\alpha_i$.

\par

Let $\Ind(\QC)$ be the set of isomorphism classes of indecomposables
of $\Mod \QC$.

\begin{conjecture}\label{conjecture}
  Let $C=\{u_1,\ldots,u_n\}$ be a cluster in a cluster algebra of
  simply-laced type and rank $n$. Let $V$ be the set of all cluster
  variables for this cluster algebra. There exists a bijection $b$ :
  $\Ind(\QC)\rightarrow V\backslash C$, $\alpha\mapsto w_\alpha$ such
  that $w_\alpha=\frac{P(u_1,\ldots,u_n)}{\prod_i u_i^{n_i}}$, where
  $P$ is a polynomial prime to $u_i$ for all $i$ and where
  $n_i=n_i(\alpha)$ is the multiplicity of the simple module
  $\alpha_i$ in the module $\alpha$.
\end{conjecture}

\begin{remark}
  Through Gabriel celebrated Theorem relating indecomposables and
  positive roots, this conjecture generalizes the Theorem of Fomin and
  Zelevinsky (\cite[Theorem 1.9]{cluster2}) which corresponds to the
  case of the alternating quiver.
\end{remark}

The main aim of the present article is to prove Conjecture
\ref{conjecture} in the case of cluster algebras of type $A_n$. This
will be done using the geometric realization in terms of
triangulations given by Fomin and Zelevinsky \cite[\S 3.5]{Ysystems}.

\section{Equivalence of categories}\label{sectionequiv}

\subsection{Triangulations and diagonals}\label{section2.1}

\begin{figure}
  \begin{center}
   \scalebox{0.75} {\includegraphics{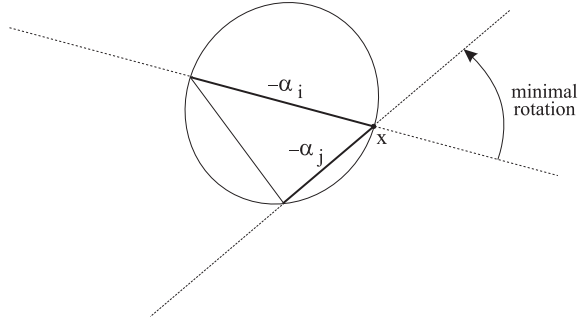}}
    \caption{$ -\alpha_i < -\alpha_j $ }
    \label{relationfig}
  \end{center}%
\end{figure}

Let us fix a nonnegative integer $n$ and a triangulation $T$ of a
regular polygon with $n+3$ vertices. The diagonals of this polygon
will be called roots and designed by Greek letters. Let us call
\textit{negative} the roots belonging to $T$ and \textit{positive} the
other roots. Let $\Phi_+$ be the set of positive roots with respect to
$T$. Let $I$ be the set of negative roots. By convention, the negative
root corresponding to $i\in I$ will be called $-\alpha_i$. If two
negative roots $-\za_i,-\za_j$ bound the same triangle in $T$, we
define a relation $<$ as follows: Denote by $x$ the common vertex of $
-\alpha_i, -\alpha_j $, then $ -\alpha_i< -\alpha_j $ if the rotation
with minimal angle around $x$ that sends the line through $-\za_i$ to
the line through $-\za_j$ is in positive trigonometric direction (see
Figure \ref{relationfig}).
The \textit{support} $\Supp \alpha \subseteq I $ of a positive root
$\alpha$ is the set of negative roots which cross $\alpha$. Note that
a positive root $\alpha$ is determined by its support. Indeed it is
possible to recover the vertices of a positive diagonal from the
sequence of crossed negative diagonals. A positive root $\alpha$ is
related to a positive root $\alpha'$ by a \textit{pivoting elementary
  move} if the associated diagonals share a vertex on the border (the
pivot), the other vertices of $\alpha$ and $\alpha'$ are the vertices
of a border edge of the polygon and the rotation around the pivot is
positive (for the trigonometric direction) from $\alpha$ to $\alpha'$.
Let $P_v$ denote the pivoting elementary move with pivot $v$. A
\textit{pivoting path} from a positive root $\alpha$ to a positive
root $\alpha'$ is a sequence of pivoting elementary moves starting at
$\alpha$ and ending at $\alpha'$.

\subsection{Categories of diagonals}\label{diagcat}

One can define a combinatorial $\CC$-linear additive category $\CatT$
as follows. The objects are positive integral linear combinations of
positive roots. By additivity, it is enough to define morphisms
between positive roots. The space of morphisms from a positive root
$\alpha$ to a positive root $\alpha'$ is a quotient of the vector
space over $\CC$ spanned by pivoting paths from $\alpha$ to $\alpha'$.

\begin{figure}
  \begin{center}
    \scalebox{0.75}{\includegraphics{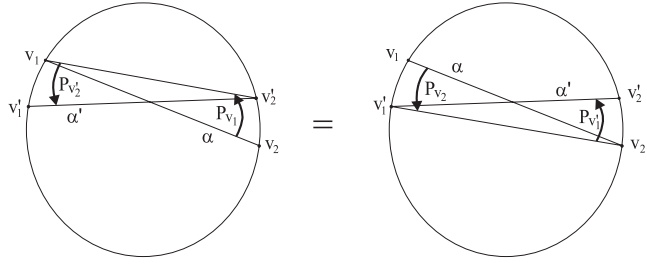}}
    \caption{mesh relation  $P_{v_2'}P_{v_1}=P_{v_1'}P_{v_2}$}
    \label{meshfig}
  \end{center}
\end{figure}

The subspace which defines the quotient is spanned by the so-called
\textit{mesh relations} (see Figure \ref{meshfig}). For any couple
$\za,\za'$ of positive roots such that $\za$ is related to $\za'$ by
two consecutive pivoting elementary moves with distinct pivots, we
define the mesh relation $P_{v_2'}P_{v_1}=P_{v_1'}P_{v_2}$, where
$v_1,v_2$ (respectively $v_1',v_2'$) are the vertices of $\za$
(respectively $\za'$) such that $P_{v_1'}P_{v_2}(\za)=\za'$. That is,
any two consecutive pivoting elementary moves using different pivots
can in some sense be ``exchanged''.

In these relations, negative roots or border edges are allowed, with
the following conventions.
\begin{itemize}
  \label{convention_mesh}
\item[(i)] If one of the intermediate edges is a border edge, the
  corresponding term in the mesh relation is replaced by
  zero.
\item[(ii)] If one of the intermediate edges is a negative root, the
  corresponding term in the mesh relation is replaced by
  zero.
\end{itemize}

More generally, a mesh relation is an equality between two pivoting
paths which differ only in two consecutive pivoting elementary moves
by such a change.

We can now define the set of morphisms from a positive root $\alpha$
to a positive root $\alpha'$ to be the quotient of the vector space
over $\CC$ spanned by pivoting paths from $\alpha$ to $\alpha'$ by the
subspace generated by mesh relations.

Therefore, the image of a pivoting path in the space of morphisms is
either the zero morphism or does only depend on the class of the
pivoting path modulo the equivalence relation on the set of pivoting
path generated by the mesh relations with no vanishing terms.

\subsection{Graphs and trees}\label{tree}

\begin{figure}
  \begin{center}
     \scalebox{0.5}{\includegraphics{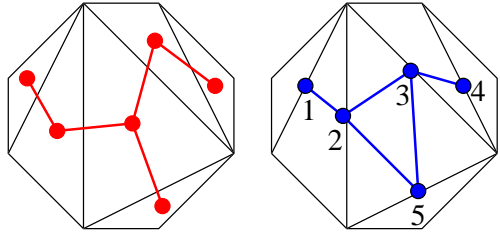}}
     \caption{Tree $t_T$ and Graph $\QT$ defined by a Triangulation
       $T$}
    \label{leaffig}
  \end{center}%
\end{figure}

Let $T$ be a triangulation. Then one can define a planar tree $t_T$ as
follows. Its vertices are the triangles of $T$ and its edges are
between adjacent triangles (see left part of Figure \ref{leaffig}).
Vertices of $t_T$ have valence $1$, $2$ or $3$. It is clear that there
is always at least one vertex of valence $1$.

From $T$, one can also define a graph $\QT$ as follows. The vertices
of $\QT$ are the inner edges of $T$ and are related by an edge if they
bound the same triangle (see right part of Figure \ref{leaffig}).

In fact, it is possible to define the graph $\QT$ starting from the
planar tree $t_T$. Vertices of $\QT$ are the edges of $t_T$. Two
vertices of $\QT$ are related by an edge in $\QT$ if the corresponding
edges of $t_T$ share a vertex in $t_T$. The equivalence with the
previous definition is obvious.

A \textit{leaf} is an edge $e$ of $t_T$ such that at least one of its
vertices has valence $1$. As there is always a vertex of valence $1$
in $t_T$, there always exists a leaf.

Recall that the mutation of a triangulation at one of its diagonal is
the unique triangulation which can be obtained by replacing this
diagonal with another one.

\subsection{Quivers}

Let $T$ be a triangulation. Let us define a quiver $\QT$ with
underlying graph the graph $\QT$ defined in the previous section.
Recall that its vertices are in bijection with $I$. Put a point at the
middle of each negative root in $T$ and draw an edge between points in
two negative roots which bound the same triangle. Two such negative
roots $-\za_i,-\za_j$ are comparable in the relation $<$ defined in
\S\ref{section2.1} and the edge is oriented $i\ot j$ if
$-\za_i<-\za_j$. From this description, it follows that all triangles
in $\QT$ are oriented. 

\begin{lemma}
  The mutation of quivers, as defined in \S\ref{sectionQCT},
  corresponds to the mutation of triangulations desdribed above.
\end{lemma}

\begin{proof} Left to the reader.
\end{proof}

One can define a $\CC$-linear abelian category
$\Mod\QT$ as follows. This is the category of modules over the quiver
$\QT$ with the following relations, called \textit{triangle
  relations}.
\begin{equation} \label{trianglerelations}
  \textup{In any triangle, the composition of two successive maps is zero.}
\end{equation}%

These relations are exactly the relations prescribed by Definition
\ref{conject-relation}.

\begin{lemma}\label{lemma-supp}
  The support $\Supp \alpha$ of a positive root $\alpha$ is connected
  as a subset of the quiver $\QT$.
\end{lemma}
\begin{proof}
  Let $-\za_{i},-\za_{j}$ be two distinct diagonals in $ \Supp\za$.
  We will show that there is an unoriented path from $i$ to $j$ in
  $\QT\cap \Supp \za$. The diagonals $-\za_i$ and $-\za_j$ cut the
  polygon into three parts. Denote by $R_{ij} $ the part that contains
  both $-\za_i$ and $-\za_j$. We proceed by induction on $m$ the
  number of negative roots in $R_{ij}$. If $m=1$ then $-\za_i=-\za_j$
  and there is nothing to prove. Let us assume that $m>1$. Let
  $\Delta$ be the unique triangle in $T$ that contains $-\za_i$ and
  lies in $R_{ij}$. Since $\za$ crosses both $-\za_i$ and $-\za_j$, it
  has to cross exactly one of the two sides different from $-\za_i$ in
  $\Delta$. This side cannot be a border edge of the polygon, hence it
  has to be a negative root, call it $-\za_k$. Thus there is an edge
  between $-\za_i$ and $-\za_k$ in $\QT$ and $i,k \in \Supp \za$. We
  may suppose by induction that there is an unoriented path in
  $\Supp\za$ from $k$ to $j$ and we are done.
\end{proof}

\subsection{Functor $\Theta$}

Let us define a $\CC$-linear additive functor $\eqv$ from $\CatT$ to
$\Mod\QT$. On objects, it is sufficient by additivity to define $\eqv$
on positive roots. The image of the positive root $\alpha$ is the
module $(M^\za,f^\za)=(M^\za_i,f^\za_{ij})$ defined by
\begin{equation*}
  M^{\alpha}_i=\left\{
    \begin{array}{ll}
      \CC &\textup{if }i \in \Supp \alpha,\\
      0 &\textup{otherwise}
    \end{array}\right.
  \quad\textup{and}\quad
  f^\za_{ij}=\left\{
    \begin{array}{ll}
      \idc &\textup{if } M_i^\za=\CC=M_j^\za,\\
      0 &\textup{otherwise.}
    \end{array}\right.
\end{equation*}
This is indeed an object in $\Mod\QT$ because a positive root $\za$
can only cross two sides of a triangle in $T$, which implies that in
each triangle in $\QT$ there is at most one arrow $i\to j$ such that
$f_{ij}^\za\ne 0$ and hence the triangle relations
(\ref{trianglerelations}) hold. Now we define the functor $\eqv$ on
morphisms. By additivity, it is sufficient to define the functor on
morphisms from a positive root to a positive root. Our strategy is to
define first the functor on pivoting elementary moves, then check that
the mesh relations hold.
For any pivoting elementary move $P\,:\,\alpha \to \alpha'$, define the
morphism $\eqv(P)$ from $(M^{\za},f^{\za})$ to $(M^{\za'},f^{\za'})$
to be $\idc$ whenever possible and $0$ else. Let us now check that
this is indeed a morphism in $\Mod\QT$. For a given arrow $j\to i$ in
$\QT$, we have to check the commutativity of the following diagram:
\begin{center}
  \begin{diagram}
    M_j^\za &\rTo^{f_{ji}^\za}& M_i^\za\\
    \dTo^{\eqv(P)_j} & & \dTo_{\eqv(P)_i}\\
    M_j^{\za'}&\rTo^{f_{ji}^{\za'}}&M_i^{\za'}
  \end{diagram}
\end{center}
This is obvious if $M_j^\za=0$ or $M_i^{\za'}=0$ and also if both
$M_i^\za$ and $M_j^{\za'}$ are $0$. Suppose $M_j^{\za}\ne 0$ and
$M_i^{\za'}\ne 0$. If $M_i^{\za}\ne 0$ and $M_j^{\za'}\ne 0$ then all
four maps are $\idc$ and the diagram commutes. The only remaining case
is if exactly one of $M_i^{\za}$, $M_j^{\za'}$ is not zero. We will
show that this cannot happen. Suppose that $M_j^{\za'}= 0$ and
$M_i^{\za}\ne 0$, that is $j,i\in \Supp \za,\ i\in \Supp \za'$ and $
j\notin\Supp\za'$. Since $\za\to \za'$ is a pivoting elementary move
we get that $-\za_{j} $ crosses $\za$, that $-\za_{j}$ and $\za'$ have
a common point on the boundary of the polygon and that $-\za_{i}$
crosses $\za$ and $\za'$. This implies that $-\za_j < -\za_{i}$ and
that contradicts the orientation $j\to i$ in the quiver $\QT$. The
other case can be excluded by a similar argument. To show that the
functor is well defined, it only remains to check the mesh relations.
Let $\za\stackrel{P^1}{\longrightarrow}\zb$,
$\zb\stackrel{P^2}{\longrightarrow}\zg$,
$\za\stackrel{P^3}{\longrightarrow}\zb'$,
$\zb'\stackrel{P^4}{\longrightarrow}\zg$ be pivoting elementary moves
with $\za,\zb,\zg$ positive roots and $\zb\ne\zb'$. Note that we can
exclude the case where $\zb$ and $\zb'$ are both negative roots or
border edges because in this case either $\za$ or $\zg$ has to be
negative too, since $T$ is a triangulation. Suppose first that $\zb'$
is positive. One has to check the commutativity of the diagram
\begin{center}
  \begin{diagram}
    M_i^\za&\rTo^{\eqv(P^1)_i}& M_i^\zb\\
    \dTo^{\eqv(P^3)_i} &&    \dTo_{\eqv(P^2)_i}\\
    M_i^{\zb'}&\rTo^{\eqv(P^4)_i}& M_i^\zg
  \end{diagram}%
\end{center}%
for all $i$. The only non trivial case is when $i\in
\Supp\za\cap\Supp\zg$. In this case, we also have $i\in
\Supp\zb\cap\Supp\zb'$ because any diagonal crossing both $\za$ and
$\zg$ must also cross $\zb$ and $\zb'$. Thus all maps are $\idc$ and
the diagram commutes. Suppose now that $\zb'$ is negative or a border
edge. We have to show that the composition $ M_i^\za
\stackrel{\eqv(P^1)_i}{\longrightarrow} M_i^\zb
\stackrel{\eqv(P^2)_i}{\longrightarrow} M_i^\zg $ is zero for all $i$.
But in this case no negative root can cross both $\alpha$ and $\zg$.
So $\Supp \alpha \cap \Supp \zg$ is empty, therefore the composition
is zero. Hence the mesh relations hold, with the conventions made in
its definition.
\begin{lemma}
  \label{lemma-homct}
  The vector space $\Hom_{\CatT}(\alpha,\alpha')$ is not zero if and
  only if there exists $i \in \Supp \alpha \cap \Supp \alpha'$ such
  that the relative positions of $\alpha$, $\alpha'$ and $-\alpha_i$
  are as in Figure \ref{homfig}. That is, let $v_1,v_2$ be the
  endpoints of $-\za_i$ and $u_1,u_2$ (respectively $u'_1,u'_2$) be
  the endpoints of $\za$ (respectively $\za'$). Then ordering the
  vertices of the polygon in the positive trigonometric direction
  starting at $v_1$, we have $v_1<u_1\le u'_1<v_2<u_2\le u_2'$. In
  this case, $\Hom_{\CatT}(\alpha,\alpha')$ is of dimension one.
%
  \begin{figure}[ht]
    \begin{center}
      \scalebox{1}{\includegraphics{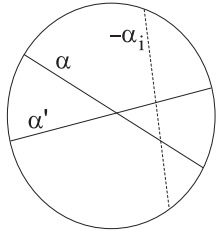}}
      \caption{Relative position}\label{homfig}
    \end{center}
  \end{figure}
\end{lemma}
\begin{proof}
  Suppose that $\Hom_{\CatT}(\alpha,\alpha')$ is not zero. Let $P\ne 0$
  be a sequence
  $\za=\za^0\stackrel{P_1}{\longrightarrow}\za^1\stackrel{P_2}
  {\longrightarrow}\ldots \stackrel{P_m}{\longrightarrow}\za^m=\za'$
  where the $\za^i$ are positive roots and the $P_i$ are pivoting
  elementary moves. Because of the mesh relations, we may suppose that
  the first $k$ moves $P_1,\ldots,P_k$ have as pivot one of the
  vertices of $\za$ and the last $(m-k)$ moves $P_{k+1},\ldots, P_m$
  have as pivot the vertex which is the intersection of $\za'$ and
  $\za^k$. Denote by $V_1$ (respectively $V_2$) the set of vertices of
  $\za^1, \ldots, \za^k$ (respectively $\za^{k+1},\ldots,\za^m$) other
  than the pivot of $P_1$ (respectively $P_{k+1}$). Since $P\ne 0$
  this implies that $\za$ and $\za'$ intersect and that all diagonals
  with one vertex in $V_1$ and the other in $V_2$ are positive. Thus
  the vertices of $\za$ and $\za'$ form a quadrilateral in the polygon
  without any diagonals of the triangulation $T$ crossing it from
  $V_1$ to $V_2$. Because $T$ is a triangulation we must have a
  $-\za_{i}\in T$ crossing this quadrilateral in the other direction
  and we get the situation in the diagram. On the other hand, in the
  situation of the diagram it is clear that there is a non-zero
  morphism from $\za$ to $\za'$. Finally, the dimension of
  $\Hom_{\CatT}(\alpha,\alpha')$ is at most one, since any two non-zero
  pivoting paths from $\alpha$ to $\alpha'$ are in the same class.
\end{proof}
\begin{lemma}\label{lemma-supp2}
  Let $\za,\za'$ be positive roots, then $\Supp \alpha \cap \Supp
  \alpha'$ is connected.
\end{lemma}%
\begin{proof}
  Suppose the contrary. Write $S=\Supp \alpha$ and $S'=\Supp \alpha'$
  for short. Let $i,k \in S\cap S'$ be two vertices that belong to
  different connected components of $S\cap S'$. Since $S$ and $S'$ are
  connected (Lemma \ref{lemma-supp}) we may choose two minimal paths
  $p: i={i_1},{i_2},\ldots,{i_p}={k}$ in $S$ and $p':
  i={j_1},{j_2},\ldots,{j_p}={k}$ in $S'$. Let $m$ be the smallest
  integer such that ${i_{m+1}}\ne{j_{m+1}} $. In the triangulation
  $T$, each of the diagonals
  $-\za_{i_{m-1}},-\za_{i_{m+1}},-\za_{j_{m+1}}$ has a vertex in
  common with $-\za_{i_m}$. Since a positive root can only cross two
  sides of a triangle, we get that $-\za_{i_m},-\za_{i_{m+1}}$ and
  $-\za_{j_{m+1}}$ form a triangle $\Delta$ in $T$. Moreover
  ${i_{m+1}}\in S\setminus S'$ and ${j_{m+1}}\in S'\setminus S$. Now
  cutting out the triangle $\Delta$ divides the polygon into three
  parts: $R_{i_m},R_{i_{m+1}}$ and $R_{j_{m+1}}$ such that $R_l$
  contains $-\za_l$, $l=i_m,i_{m+1},j_{m+1}$. Clearly all $-\za_{i_l},
  \, l\ge m+1$ lie in $R_{i_{m+1}}$ and all $-\za_{j_l}, \, l\ge m+1$
  lie in $R_{j_{m+1}}$. But this contradicts the fact
  $-\za_{j_q}=-\za_k=-\za_{i_p}$ and we have shown that $S\cap S'$ is
  connected.
\end{proof}%
\begin{lemma}
  \label{lemma-homqt}
  The vector space
  $\Hom_{\Mod\QT}((M^{\za},f^{\za}),(M^{\za'},f^{\za'}))$ is not zero
  if and only if the following conditions hold. Let $S=\Supp \alpha$
  and $S'=\Supp \alpha'$ for short.
  \begin{itemize}
  \item[\textup{(i)}] $S \cap S'$ is not empty,
  \item[\textup{(ii)}] There is no arrow from $S \setminus S'$ to $S \cap S'$
    in $\QT$,
  \item[\textup{(iii)}] There is no arrow from $ S \cap S'$ to $S' \setminus S
    $ in $\QT$.
  \end{itemize}
  In this case, $\Hom_{\Mod\QT}((M^{\za},f^{\za}),(M^{\za'},f^{\za'}))$
  is of dimension one.
\end{lemma}
\begin{proof}
  Let $P$ be a non-zero element of
  $\Hom_{\Mod\QT}((M^{\za},f^{\za}),(M^{\za'},f^{\za'}))$. Then
  condition (i) is clearly true. Let us show conditions (ii) and
  (iii). Suppose that condition (ii) is not true, thus there is an
  arrow $i\to j$ in $\QT$ with $i\in S\setminus S'$, $j\in S\cap S'$
  and such that the following diagram commutes.
  \begin{center}
    \begin{diagram}
      M_i^\za&\rTo^{\idc}& M_j^\za\\
      \dTo^{P_i} &&    \dTo_{P_j}\\
      0&\rTo& M_j^{\za'}
    \end{diagram}%
  \end{center}%
  Thus $P_i$ and $P_j$ are both zero. Now let $k$ be any vertex in
  $S\cap S'$, We will show that $P_k$ is zero. By Lemma
  \ref{lemma-supp2}, there is an unoriented path $k=k_0
  \frac{\quad}{\quad} k_1 \frac{\quad}{\quad} \cdots$
  $\frac{\quad}{\quad} k_m=i$ in $\QT$ such that each $k_i\in S\cap
  S'$. We proceed by induction on $m$. The case $m=0$ is done above,
  suppose $m>0.$ By induction $P_{k_1}$ is zero and the commutativity
  of the diagram
  \begin{center}
    \begin{diagram}
      M_k^\za&\rLine^\idc& M_{k_1}^\za\\
      \dTo^{P_k} &&    \dTo_{0}\\
      M_k^{\za'}&\rLine^\idc& M_{k_1}^{\za'}
    \end{diagram}%
  \end{center}%
  implies that $P_k$ is zero for both possibilities of orientation of
  $k\frac{\quad}{\quad} k_1$ in $\QT$. By contradiction, this shows
  (ii). Condition (iii) is proved by a similar argument.
  
  In order to show the converse statement, let $\za, \za'$ be such
  that (i),(ii) and (iii) hold. Define
  $P\in\Hom_{\Mod\QT}((M^{\za},f^{\za}),(M^{\za'},f^{\za'}))$ by
  $P_i=\idc$ whenever $i\in S\cap S'$ and $P_i=0$ otherwise. Then (i)
  implies that $P$ is non-zero. We only have to check that $P $ is a
  morphism of quiver modules, \textit{i.e.} that the diagram
  \begin{center}
    \begin{diagram}
      M_i^\za&\rTo^{f_{ij}^\za}& M_{j}^\za\\
      \dTo^{P_i} &&    \dTo_{P_j}\\
      M_i^{\za'}&\rTo^{f_{ij}^{\za'}}& M_{j}^{\za'}
    \end{diagram}%
  \end{center}%
  commutes for all $i\to j$ in $\QT$. But this is true because of
  conditions (ii) and (iii). Finally, the dimension of
  $\Hom_{\Mod\QT}((M^{\za},f^{\za}),(M^{\za'},f^{\za'}))$ is at most
  one, since all vector spaces $M^{\za}_i,M^{\za'}_i$ are of dimension
  zero or one, and by connexity of the intersection of the supports
  (Lemma \ref{lemma-supp2}).
\end{proof}

\begin{lemma}\label{lemma-homeqv}
  The conditions in Lemma \ref{lemma-homct} and in Lemma
  \ref{lemma-homqt} are equivalent.
\end{lemma}
\begin{proof}
  Suppose $\za,\za'$ are as in Lemma \ref{lemma-homct}. Then $i\in
  S\cap S'$ which implies (i). Suppose there is an arrow $j \to k$ in
  $\QT$ such that $j\in S\setminus S'$ and $k \in S\cap S'$. Then
  $-\za_{k}$ crosses both $\za$ and $\za'$ while $-\za_{j}$ crosses
  only $\za$. Since $j\to k$, we know that $-\za_{j}$ and $-\za_{k}$
  are two sides of a triangle in $T$ and that $-\za_{k} < -\za_{j}$.
  This is impossible because of the way that $-\za_{i}$ intersects
  $\za$ and $\za'$. We have shown that condition (ii) holds;
  condition (iii) can be shown similarly. This proves one direction of
  the Lemma.
  
  Suppose now that $\za,\za'$ satisfy the conditions (i),(ii) and
  (iii) of Lemma \ref{lemma-homqt}. By (i), there exists $-\za_i$ in
  $S \cap S'$. Let $v_1,v_2$ be the endpoints of $-\za_i$. Consider
  the two parts of the polygon $R_l$ and $R_r$ delimited by
  $-\alpha_{i}$. Each of them contains exactly one vertex of $\alpha$
  and exactly one vertex of $\alpha'$. Consider the positive roots
  $\za,\za'$ as a paths running from $R_r$ to $R_l$. The sequence of
  negative roots given by the successive intersections of the path
  $\za$ (respectively $\za'$) with elements of $T$ yields an ordering
  of $\Supp\za$ (respectively $\Supp\za'$). Let
  $S_l=\{-\za_i=-\za_{i_1},-\za_{i_2},\ldots,-\za_{i_p}\}$
  (respectively
  $S'_l=\{-\za_i=-\za_{j_1},-\za_{j_2},\ldots,-\za_{j_q}\}$) be the
  set of negative roots in $R_l$ crossing $\za$ (respectively $\za'$)
  in that order. Let $m$ be the greatest integer such that
  $-\za_{i_m}=-\za_{j_m}$. We will distinguish four cases.
  \begin{itemize}
  \item[1.] $m=p=q$, then on the boundary of $R_l$, going from the
    endpoints of $-\za_i$ in positive direction, we meet $\za$ and
    $\za'$ at the same time.
  \item[2.] $m=p<q$, then $-\za_{j_{m+1}}$ and $-\za_{j_m}$ bound the
    same triangle in $T$. The corresponding edge in $\QT$ is oriented
    $j_{m+1}\to j_m$ by (iii). This implies that going from the
    endpoints of $-\za_i$ in positive direction on the boundary of
    $R_l$, we meet $\za$ first and then $\za'$.
  \item[3.] $m=q<p$, then $-\za_{i_{m+1}}$ and $-\za_{i_m}$ bound the
    same triangle in $T$. The corresponding edge in $\QT$ is oriented
    $i_{m+1}\ot i_m$ by (ii). This implies again that going from the
    endpoints of $-\za_i$ in positive direction on the boundary of
    $R_l$, we meet $\za$ first and then $\za'$.
  \item[4.] $m<p$ and $m<q$, then $-\za_{i_{m+1}}$, $-\za_{i_m}$ and
    $-\za_{j_{m+1}}$ are three different diagonals that bound the same
    triangle in $T$. The corresponding edges in $\QT$ are oriented
    $i_{m+1}\ot i_m$ by (ii) and $j_{m+1}\to j_m$ by (iii). This
    implies once more that going from the endpoints of $-\za_i$ in
    positive direction on the boundary of $R_l$, we meet $\za$ first
    and then $\za'$.
  \end{itemize}
  By symmetry, we obtain the same results in the other part $R_r$.
  This implies that the relative position of $\za,\za'$ and $\za_i$ is
  exactly the one described in Lemma \ref{lemma-homct}.
\end{proof}

\begin{proposition}
  The functor $\eqv$ is fully faithful.
\end{proposition}
\begin{proof}
  Using Lemmas \ref{lemma-homct}, \ref{lemma-homqt} and
  \ref{lemma-homeqv}, it only remains to show that the image of a
  non-zero morphism is a non-zero morphism. It is sufficient to show
  this for all non-zero morphisms between positive roots. Let $P\in
  \Hom(\za,\za')$ be such a morphism. Then $P$ is given by a sequence
  of pivoting elementary moves $\ \za=\za^1\to \ldots\to\za^m=\za'$.
  This sequence being a non-zero morphism implies that there exists a
  negative root $-\za_i$ crossing all the $\za^k$, $k=1,\ldots ,m$, by
  Lemma \ref{lemma-homct}. By definition, $\eqv(P)_i $ is $\idc$,
  hence non-zero.
\end{proof}


\begin{remark}
  \label{leafremark}
  If $i$ is a leaf (see \S\ref{tree}) and $\za$ a positive root, then
  $i\in \Supp \za$ if and only if one endpoint of $\za$ is the vertex
  $x$ of the polygon that is cut off by $-\za_i$. Therefore there
  exists one positive root $\za^{Pr_i}$ such that the set of all
  positive roots with $i$ in their supports is equal to the set
  \begin{equation*}
    \{\za^{Pr_i},P_x(\za^{Pr_i}),P_x(P_x(\za^{Pr_i})),\dots,
  P_x^{n-1}(\za^{Pr_i})\},
  \end{equation*}
  where $P_x$ is the pivoting elementary move with pivot $x$.
\end{remark}

\begin{theorem}\label{equivalence}
  The functor $\eqv$ gives an equivalence of categories from $\CatT$ to
  $\Mod\QT$.
\end{theorem}
\begin{proof}
  It only remains to show that the functor $\eqv$ is essentially
  surjective, \textit{i.e.} that each indecomposable module in
  $\Mod\QT$ is the image of a positive root under $\eqv$. In fact, we
  will characterize the indecomposable modules of $\Mod\QT$ with the
  help of the Auslander-Reiten theory, and this will enable us to
  conclude by proving that there are $\frac{n(n+1)}{2}$ indecomposable
  $\QT$-modules.

  In the following, we refer to \cite{gabriel}, see also \cite{AR} for
  definitions, notation and results in representation theory of finite
  dimensional algebras. In this proof, we use the following notations
  in $\Mod\QT$: $P_i$ (respectively $I_i$) is the $i-$th projective
  (respectively injective) indecomposable module, with the convention
  that $P_{\emptyset}=I_{\emptyset}=0$. Hence $(P_i)_l=\CC$ if there
  is an oriented path in $\QT$ modulo the triangle relations from $i$
  to $l$ and $(P_i)_l=0$ otherwise. Similarly, $(I_i)_l=\CC$ if there
  is an oriented path in $\QT$ modulo the triangle relations from $l$
  to $i$ and $(I_i)_l=0$ otherwise. The maps of $P_i$ and of $I_i$ are
  $\idc$ whenever possible and zero otherwise. In particular, these
  modules are multiplicity free. Fix a triangulation $T$.  In the
  sequel, we set $Q=\QT$ when no confusion occurs and we denote by
  $Q_0$ its set of vertices. Given a subset of vertices $S$ of $Q_0$,
  a {\it full subquiver} of $Q$ with vertices $S$ will be the set $S$
  together with the set of all arrows (with relations) of $Q$ joining
  vertices of $S$. We say that a $\QT$-module $M$ is {\it of type }$A$
  if the full subquiver of $\QT$ on the support of $M$ is of type
  $A_k$ for some $k\in{\NN}$.

\begin{lemma}
  Let $M$ be an indecomposable $Q$-module of type $A$ and let $N$ be any
  indecomposable $Q$-module. If $\Hom_Q(N,M)$ or $\Hom_Q(M,N)$
  contains an irreducible morphism, then $N$ is of type $A$.
\end{lemma}
\begin{proof}
  The proof is based on the construction of irreducible morphisms via
  the Nakayama functor \cite[\S 4.4]{gabriel}. The dual functor gives
  an (anti)-equivalence between $\Mod Q$ and $\Mod Q^{opp}=\Mod
  Q_{T^*}$, where $T^*$ is the mirror triangulation. Using this
  (contravariant) functor, we can easily reduce the proof to the case
  where $\Hom_Q(N,M)$ contains an irreducible.
  
  Suppose that the support of $M$ is given by the set
  $Q_0':=\{1,\ldots,m\}$. Let $Q'$ be the full subquiver of $Q$ given
  by $\Supp M$. By assumption, $Q'$ is of type $A_m$ with extremal
  vertices 1 and $m$ and we can suppose that the edges link $i$ with
  $i\pm 1$. Remark that, as $M$ is indecomposable of type $A$, $M$ is
  multiplicity free.
  
  \begin{figure}[ht]
    \begin{center}
      \scalebox{0.75}{\includegraphics{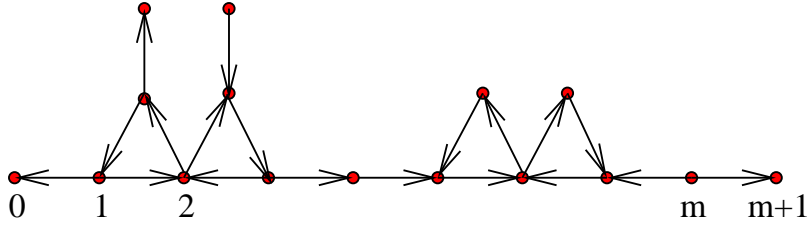}}
      \caption{Quiver $Q$ and subquiver $Q'$ of type $A$}\label{nakayama}
    \end{center}
  \end{figure}

  In order to go further, we need some more precise description of the
  quiver $Q'$ inside $Q$. The reader may like to follow the argument
  on the example provided by Figure \ref{nakayama}. Suppose that the module 
$M$ is not simple. Note that there
  exists at most one vertex which will be denoted by $0$, resp. $m+1$,
  in $Q_0\backslash Q_0'$ such that $1\rightarrow 0$, resp.
  $m\rightarrow m+1$, and such that there exist no other edges between
  $Q_0'$ and $0$, resp. $m+1$. Note also that for all $k$ in $Q_0'$,
  there exist at most two vertices $i(k)^{\pm}$ such that
  $i(k)^{\pm}\not\in Q_0'$ and $i(k)^{\pm}\rightarrow k$ is an arrow
  of $Q$. We can define $i(k)^+$, resp. $i(k)^-$, to be the vertex
  such that $k+1\rightarrow i(k)^+$, resp. $k-1\rightarrow i(k)^-$. By
  convention, if the vertices $i(k)^+$, $i(k)^-$, $0$, $m+1$ do not
  exist, we define the corresponding symbol to be the empty set.

  Let $S_{0}$, resp. $S_{m+1}$, be the support of the injective
  module associated to $0$, resp. $m+1$, in the full subquiver of $Q$
  with set of vertices $Q\backslash{1}$, resp. $Q\backslash{m}$. It is
  clear from the tree structure of $Q$, see \S\ref{tree}, that the set
  $Q_0'':=S_0\cup Q_0'\cup S_{m+1}$ is the set of vertices of a full
  subquiver $Q''$ of $Q$ of type $A$.

  For each source $k$ or for $k=1,m$, let $S_k^{\pm}$ be the support
  of the injective module associated to $i(k)^{\pm}$ in the full
  subquiver of $Q$ with set of vertices $Q_0\backslash{Q'_0}$.

  Note that if $M$ is projective, then the module $N$ is a direct summand of 
the radical of $M$ and the lemma is true in this case. Suppose now that $M$ is not projective. 
  Then, we can calculate the Auslander-Reiten translate ${\mathcal
    A}(M)$ via the Nakayama functor. We have a minimal projective
  presentation of $M$: $P^1\rightarrow P^0\rightarrow M\rightarrow 0$
  by setting $P^0:=\oplus_{i\in R}P_i$, where $i$ runs over the set
  $R$ of sources of $Q'$, and $P^1:=\oplus_{j\in R'}P_j$, where $j$
  runs over the set $R'$ given by the union of $\{0,m+1\}$ with the
  sinks of $Q_0'\backslash\{1,m\}$.
  
  By \cite[Remark 3.1]{gabriel}, with the help of the Nakayama
  functor, we obtain a minimal injective representation of the
  Auslander-Reiten translate ${\mathcal A}(M)$ of $M$: $0\rightarrow
  {\mathcal A}(M)\rightarrow I^1\rightarrow I^0$, with
  $I^0:=\oplus_{i\in R}I_i$, and $I^1:=\oplus_{j\in R'}I_j$. For $k$
  in $R'$, let $k^+$, resp. $k^-$, be the source of $Q'$ succeeding,
  resp. preceding, $k$, with the convention that $k^+=m$, resp.
  $k^-=1$, if there is no such source.

  Then, for each sink $k$ of $Q_0'$, the support of $I_k$ is given by
  $\Supp (I_k)=S_{k^-}^-\cup S_{k^+}^+\cup(\Supp (I_k)\cap Q_0')$. The
  support of $I_0$, resp. $I_{m+1}$, is $S_{l}^+\cup(\Supp (I_0)\cap
  Q_0'')$, resp. $S_{g}^-\cup(\Supp (I_{m+1})\cap Q_0'')$, where $l$,
  resp. $g$, is the lowest, resp. greatest, source of $Q_0'$. The
  support of $I_k$, $k\in R$ contains $S_k^{\pm}$.
  
  This implies that the support of ${\mathcal A}(M)$ is a subset of
  $Q_0''$.
  
  Now, let $X$ be the middle term in the Auslander-Reiten sequence
  $0\rightarrow {\mathcal A}(M)\rightarrow X\rightarrow M\rightarrow
  0$. We have $\Supp X\subset\Supp{\mathcal A}(M)\cup\Supp M\subset
  Q_0''$. Hence, $X$ is of type $A$. The Auslander-Reiten
  Theorem asserts that the module $N$ of the Lemma is a direct summand
  of $X$. So, we obtain the Lemma in this case. The case where $M$ is
  simple is very similar and left to the reader.
\end{proof}
We can now prove the Theorem. Let $M$ be an indecomposable $Q$-module
of type $A$. By the Lemma, the component of $M$ in the
Auslander-Reiten quiver contains only modules of type $A$, therefore
can only be finite. Hence, by \cite[Proposition 6.3]{gabriel}, every
indecomposable module is of type $A$ and in particular is multiplicity
free. So, there exists a one-to-one correspondence between
indecomposable $Q$-modules and full subquivers of $Q$ of type $A$. Let
$i$ be a leaf of $Q$ and $j$ be any vertex of $Q$. By \S\ref{tree},
there exists a unique full subquiver of type $A$ of $Q$ whose extreme
vertices are $i$ and $j$. This implies by induction that the number of
such subquivers is $\frac{n(n+1)}{2}$. Hence, there are
$\frac{n(n+1)}{2}$ indecomposable $Q$-modules as required.

\end{proof}%

\begin{corollary}
  The category $\CatT$ is abelian.
\end{corollary}

\begin{corollary}
  \label{maincor}
  There exists a bijection $\varphi$ between $\Ind(\QT)$ and the
  diagonals of the polygon not in $T$. Moreover, for $M$ in
  $\Ind(\QT)$ and any vertex $i$ of $\QT$, the multiplicity of the
  simple module $S_i$ in the module $M$ is 1 if $\varphi(M)$ crosses
  the $i^{th}$ diagonal of $T$ and 0 if not. In particular, for two
  isoclasses $M$, $M'$ in $\Ind(\QT)$, we have $M=M'$ if and only if
  $n_i(M)= n_i(M')$ for all $i$.
\end{corollary}

\subsection{The orbit category}

This subsection is not used in the sequel. We give here a description
of the category $\CatT$, using the equivalence of category proved
above. Then, we prove that the orbit category introduced by
\cite{marall} has a nice geometric realization in the $A_n$ case. Let
$r^+$, resp. $r^-$, be the elementary rotation of the polygon in the
positive, resp. negative, direction.
\begin{theorem}
  \label{summary}
  Let $T$ be a triangulation of the $n+3$ polygon, and let $\CatT$ be
  the corresponding category, then:
  \begin{itemize}
  \item[\textup{(i)}] The irreducible morphisms of $\CatT$ are direct
    sums of the generating morphisms given by pivoting elementary
    moves.
  \item[\textup{(ii)}] The mesh relations of $\CatT$ are the mesh
    relations \cite{AR} of the Auslander-Reiten quiver of $\CatT$.
  \item[\textup{(iii)}] The Auslander-Reiten translate is given on
    diagonals by $r^-$.
  \item[\textup{(iv)}] The projective indecomposable objects of
    $\CatT$ are diagonals in $r^+(T)$.
  \item[\textup{(v)}] The injective indecomposable objects of $\CatT$
    are diagonals in $r^-(T)$.
  \end{itemize} 
\end{theorem}
\begin{proof}
  (i) and then (ii) are clear by construction of the category $\CatT$.
  By (i) and (ii), extremal terms of an almost split sequence are
  given by the diagonals $\alpha$ and $\alpha'$ of Figure 2. This
  proves (iii). (iv) and (v) follows from (iii).
\end{proof}

The assertions (iv) and (v) of Theorem \ref{summary} suggest an
interpretation of the diagonals of the triangulation $T$ in terms of
category. Indeed, we will consider those diagonals, at least in the
hereditary case, as shifts of the projectives in the derived category
${\cal D}\CatT$.

In order to simplify the construction, suppose that $T$ is a
triangulation corresponding to the unioriented quiver $A_n$ with
simple projective $P_1$. The category $\CatT\simeq\Mod\QT$ is
hereditary, so, the indecomposable objects of the derived category
${\mathcal D}\Mod\QT$ are the shifts $M[m]$, $m\in \ZZ$, of the
indecomposables $M$ of $\Mod\QT$. Let $F$ be the functor of ${\mathcal
  D}\Mod\QT$ given by $M\mapsto{\mathcal A}^{-1}(M)[1]$, where
${\mathcal A}$ is the Auslander-Reiten translate in the derived
category. We define the orbit category $\overline{{\cal D}\Mod\QT}$
whose objects are objects of ${\cal D}\Mod\QT$ and morphisms are given
by $\Hom_{\overline{{\cal D}\Mod\QT}}(M,N):=\oplus_{i\in
  \ZZ}$$\Hom_{{\cal D}\Mod\QT}(M,F^i(N))$, $M$, $N\in {\mathcal
  D}\Mod\QT$. The set $\Ind(\QT)\cup\{P_i[1],\,1\leq i\leq n\}$ is the
set of indecomposable objects of $\overline{{\mathcal D}\Mod\QT}$ up
to isomorphism. Note that the category $\overline{{\mathcal
    D}\Mod\QT}$ is triangulated but not abelian in general.

We can also construct the total category $\mathcal{C}$ generated by
all the diagonals of the $(n+3)$ polygon. The construction is analogue
to the construction of $\CatT$: indecomposable objects are positive
roots and simple negative ones. The homomorphisms and the mesh
relations are defined as in \S\ref{diagcat} without the point (ii) in
the convention made there.

\begin{theorem}\label{orbit}
  The categories $\mathcal{C}$ and $\overline{{\mathcal D}\Mod\QT}$
  are equivalent.
\end{theorem}
We give here a sketch of the proof. The derived category of
representations of the unioriented quiver is well-known. The
indecomposable objects of ${\mathcal D}\Mod\QT$ can be indexed by
$\kappa$: $\ZZ \times\{1,\ldots,n\}\rightarrow\Ind{\mathcal
  D}\Mod\QT$, by the rule $\kappa(1,i)=P_i$, $\kappa(i+1,j)={\mathcal
  A}^{-1}(\kappa(i,j))$, $\kappa(i+j,n+1-j)=\kappa(i,j)[1]$. This
implies that the indecomposable objects of $\overline{{\mathcal
    D}\Mod\QT}$ can be indexed by $\ZZ \times\{1,\ldots,n\}\slash
(i,j)\equiv (i+j+1,n-j+1)$.

For all $(i,j)$ in $\ZZ\times\{1,\ldots,n\}$, we define the
quadrilateral $R_{(i,j)}$ by its vertices $(i,j)$, $(i,n)$,
$(i+j-1,1)$, $(i+j-1,n-j+1)$. Let $M$ be the indecomposable object in
$\overline{{\mathcal D}\Mod\QT}$ indexed by $(i,j)$, an let $N$ be any
indecomposable object. Then, $\Hom_{\overline{{\mathcal
      D}\Mod\QT}}(M,N)\not=0$ if and only if $N$ is indexed by a point
inside $R_{(i,j)}$. In this case, it is ${\CC}$ as a space and the
composition of morphisms is given by the multiplication. Now, let us
index the vertices of the $(n+3)$ polygon by the group $\ZZ/(n+3)\ZZ$
and let $[i,j]$ be the diagonal from $i$ to $j$, $j-i\not=1,0,-1$. By
the description above, the additive functor defined by
$[i,j]\mapsto\kappa(i,j-i-1)$ gives an equivalence of category.

\begin{remark}
  Using the equivalence above, it is easy to see that given two
  diagonals $\alpha$ and $\alpha'$, the group $\Ext^1(\alpha,\alpha')$
  is non zero if and only if $\alpha$ and $\alpha'$ cross. Hence, a
  triangulation of the polygon correspond to a maximal set of pairwise
  extension free diagonals. 
\end{remark}

\begin{remark}
  The orbit category $\overline{{\mathcal D}\Mod\QT}$ was introduced
  in \cite{marall} for all simply-laced root systems. Its construction
  was given to us by Bernhard Keller.
\end{remark}

\begin{figure}
  \begin{center}
    \scalebox{0.5}{\includegraphics{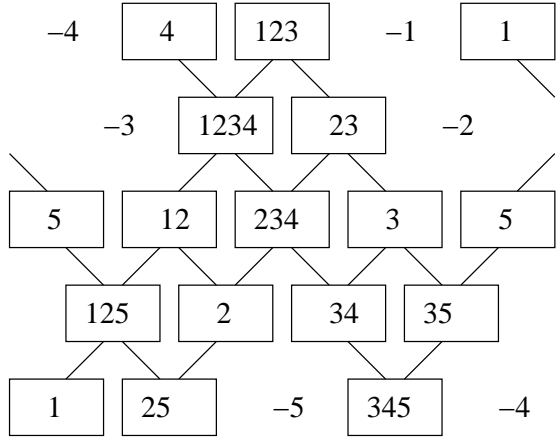}}
    \caption{Example of Auslander-Reiten Quiver}\label{ARqui}
  \end{center}
\end{figure}

An example of the Auslander-Reiten quiver is provided in Figure
\ref{ARqui}, for the Quiver with relation shown in the right part of
Figure \ref{leaffig}.

\section{Denominators of Laurent polynomials}

In this section, we will prove Conjecture \ref{conjecture} for cluster
algebras of type $A$. The existence of the bijection will be an easy
consequence of section \ref{sectionequiv} and older results of Fomin
and Zelevinsky. Most of this section is concerned with the calculation
of the exponents in the denominators of the Laurent polynomials, in
order to prove the equation of the Conjecture.

Throughout this section we will use the following setup. Let
$\Phi_{\ge -1}$ be the set of almost positive roots and $\tau_+$ and
$\tau_-$ the involutions on $\Phi_{\ge-1}$ defined by Fomin and
Zelevinsky in \cite{Ysystems}. For any $\za\in\Phi_{\ge -1}$, let
$w_\za$ be the cluster variable corresponding to $\za$ by the
bijection of Fomin and Zelevinsky \cite[Theorem 1.9]{cluster2}. Let
$C=\{u_1,\ldots,u_n\}$ be a cluster and let $\zb_1,\ldots,\zb_n$ be
the almost positive roots such that $w_{\zb_i}=u_i$. Recall the
following properties of $\tau_\pm$:

\begin{proposition}
  \label{tau} 
  \begin{itemize}
  \item[\textup{1.}] {Every $\langle \tau_-,\tau_+ \rangle$--orbit in
      $\Phi_{\ge-1}$ contains a negative simple root.}
  \item[\textup{2.}] {There is a unique function $(\,\mid\mid\,)
      :\Phi_{\ge-1}\times\Phi_{\ge-1}\to \ZZ_{\ge 0} $ such that}
    \begin{itemize}
    \item[\textup{(i)}] $(-\za_i\mid\mid\zb)= \max(n_i(\zb),0)$,
    \item[\textup{(ii)}]$(\tau_\pm\za\mid\mid \tau_\pm\zb)=(\za\mid\mid\zb).$
    \end{itemize}
    Furthermore $(\,\mid\mid\,)$ is symmetric for simply-laced root systems.
  \item[\textup{3.}]{The set
      $\tau_\pm(C)=\{w_{\tau_\pm(\zb_i)}\mid i=1,\ldots,n\}$ is a
      cluster.}
  \end{itemize}
\end{proposition}
\begin{proof}
  1. is shown in \cite[Theorem 2.6]{Ysystems}, 2. in \cite[section
  3.1]{Ysystems} and 3. in \cite[Proposition 3.5]{Ysystems}.
\end{proof}%
By the Laurent phenomenon \cite{cluster1}, we can write, for any almost
positive root $\za$,
\begin{equation}\label{eqlaurent}
  w_\za=\frac{R_{\za,C}}{\prod_{i=1}^n u_i^{[\za,\zb_i,C]}},
\end{equation}
where $R_{\za,C}$ is a polynomial in the variables $u_1,\ldots,u_n$
such that none of the $u_i$ divides $R_{\za,C}$, and
$[\za,\zb_i,C]\in\ZZ$. The following Lemma is crucial.
\begin{lemma}\label{cruciallemma} 
  For any pair of almost positive roots $\za,\zb_i$ and any pair of
  clusters $C,C'$ such that $u_i=w_{\zb_i} \in C\cap C'$, we have
  \begin{equation*}
    [\za,\zb_i,C]=[\za,\zb_i,C'].
  \end{equation*}
\end{lemma}
\begin{proof}
  It is sufficient to prove assertions (a) and (b) below.
  \begin{itemize}
  \item[(a)] All clusters containing the given cluster variable $u_i$
    are connected in the mutation graph.
  \item[(b)] In mutations which do not exchange $u_i$, the exponent of
    $u_i$ in the denominator of $w_\za$ is unchanged.
  \end{itemize}%

  \noindent
  
  Assertion (a) can either be seen as a classical statement on the
  link of a simplex in a simplicial sphere, or can be checked directly
  using the recursive properties of clusters. First the adjacency
  graph of clusters containing a fixed cluster variable is mapped to
  an isomorphic graph by the action of $\tau_+$ and $\tau_-$ because
  this action respects the compatibility function. Thus one can
  suppose that the fixed cluster variable is a negative one. But then
  the graph of clusters containing this negative simple root is
  isomorphic to a product of whole cluster adjacency graphs for
  smaller sub-root systems. Therefore it is connected as a product of
  connected graphs.

  
  Let us show (b) now. Consider the mutation $C\to C'$ that exchanges
  the cluster variables $u_j\in C $ and $u_j'\in C'$. The exchange
  relation gives $u_j=\frac{M_1+M_2}{ u_j'}$, where $M_1$ and $M_2$ are
  monomials without common divisors in the variables $C\setminus
  \{u_j\}$. From equation (\ref{eqlaurent}), we obtain by substitution
  \begin{equation*}
    w_\za=
    \frac
    {R_{\za,C}(u_1,\ldots,u_{j-1},\frac{M_1+M_2}{u'_j},u_{j+1},\ldots,u_n)}
    {\prod_{l\not=j} u_l^{[\za,\zb_l,C]}
      {(\frac{ M_1+M_2}{u'_j})}^{[\za,\zb_j,C]}}.
  \end{equation*}
  By the Laurent phenomenon, we know that $w_\za$ is a Laurent
  polynomial in the cluster variables $u_1,\ldots,u'_j,\ldots,u_n$. We
  want to prove that the exponent of $u_i$ in the denominator is still
  $[\za,\zb_i,C]$. Clearly, by properties of $M_1$ and $M_2$, this is
  true if and only if the Laurent polynomial
  $R_{\za,C}(u_1,\ldots,\frac{M_1+M_2}{ u'_j},\ldots,u_n)$ is not zero
  after evaluation at $u_i=0$. To conclude, remark that, by properties
  of the monomials $M_1$ and $M_2$ stated above, the value of this
  Laurent polynomial at $u_i=0$ is obtained by an invertible
  substitution from the value of $R_{\za,C}$ at $u_i=0$, which is
  known not to be zero.
\end{proof}

Thus $[\za,\zb,C]$ will be denoted simply $[\za,\zb]$ from now on.


The following Lemma is proved for all simply-laced root-systems.

\begin{lemma}\label{taulemma}
  Let $\alpha$,$\beta$ be two almost positive roots. Then
  \begin{equation}
    [ \alpha , \beta ]= [\tau_{\pm} \alpha , \tau_{\pm}\beta ].
  \end{equation}
\end{lemma}

\begin{proof}

  Let us consider a sequence of adjacent clusters
  \begin{equation}
    C_0 \leftrightarrow C_1 \leftrightarrow \dots \leftrightarrow C_N,
  \end{equation}
  where $\alpha \in C_0$ and $\beta \in C_N$. The exchange relations
  depend only on the matrices $B(C_0),\dots,B(C_N)$ associated to
  these clusters. As the action of $\tau$ respects the compatibility
  relation, one gets another chain of adjacent clusters
  \begin{equation}
    \tau (C_0) \leftrightarrow \tau(C_1) \leftrightarrow \dots
    \leftrightarrow \tau(C_N),
  \end{equation}
 where $\tau(\alpha) \in \tau(C_0)$ and $\tau(\beta)\in \tau(C_N)$.
 By Lemma 4.8 in \cite{cluster2}, one has
 \begin{equation}
   B_{\tau \gamma',\tau \gamma''}(\tau(C))=-B_{\gamma',\gamma''}(C),
 \end{equation}
 for any cluster $C$ and roots $\gamma',\gamma''$ in it. This minus
 sign does not change the exchange relations. From this one deduces
 that the expression of the cluster variable $u_\alpha$ in the
 variables of the cluster $C_N$ is the same as the expression of the
 cluster variable $u_{\tau(\alpha)}$ in the variables of the cluster
 $\tau(C_N)$. This proves the Lemma.
\end{proof}

\begin{lemma}
  \label{simplelemma}
  Let $-\za_i$ be a simple negative root and $\za$ an almost positive
  root. Then
  \begin{equation*}
    [\za,-\za_i]=n_i(\za).
  \end{equation*}
\end{lemma}%
\begin{proof}
  The quantity $[\za,-\za_i]$ is computed using the expression of the
  cluster variables $\za$ in the cluster made of negative roots. Then
  the bijection of Fomin and Zelevinsky between cluster variables and
  roots (\cite[Theorem 1.9]{cluster2}) implies that this is
  $(\za\mid\mid-\za_i)$. By Proposition \ref{tau}.2.(i) and symmetry
  of $(\,\mid\mid\,)$ in the simply-laced cases, the conclusion
  follows.
%
\end{proof}%
\begin{proposition}\label{prop}
  Let $\alpha$,$\beta$ be two distinct almost positive roots. Then
  \begin{equation*}
    [ \alpha , \beta ]= ( \alpha \mid\mid \beta).
  \end{equation*}
\end{proposition}
\begin{proof}
  Define a function $b\,:\, \Phi_{\ge-1}\times \Phi_{\ge-1}\to
  \ZZ_{\ge 0}$ by
  \begin{equation*}
    b(\za,\zb)=\left\{\begin{array}{ll}
  [\za,\zb] &\textup{if }\za\ne\zb,\\
  0&\textup{if }\za=\zb.
 \end{array}\right.
  \end{equation*}
  This function is well defined by Lemma \ref{cruciallemma}. Moreover
  \begin{itemize}
  \item[(1)] $b(-\za_i,\zb)   =\max(n_i(\zb),0)$ (by Lemma
    \ref{simplelemma})
  \item[(2)] $b(\tau_\pm\za,\tau_\pm\zb)=b(\za,\zb)$
    (by Lemma \ref{taulemma}).
  \end{itemize}%
  By Proposition \ref{tau}(2), the function $(\ \mid\mid\ 
  ):\Phi_{\ge-1}\times\Phi_{\ge-1}\to \ZZ_{\ge 0} $ is the unique
  function having the properties (1) and (2), thus
  $b(\za,\zb)=(\za\mid\mid\zb)$. Therefore
  $[\za,\zb]=(\za\mid\mid\zb)$ if $\za\ne\zb$.
\end{proof}
The following Theorem establishes the Conjecture \ref{conjecture} for
the type $A_n$.
\begin{theorem}
  Let $C=\{u_1,\ldots,u_n\}$ be a cluster of a cluster algebra of type
  $A_n$ and let $V$ be the set of all cluster variables of the
  algebra. Let $\QC$ be the quiver with relations associated to $C$
  and $\Ind(\QC)$ the set of isoclasses of indecomposable modules.
  Then there is a bijection
  \begin{equation*}
    \Ind(\QC) \to V\setminus C, \ \za\mapsto w_\za,
  \end{equation*}
  such that
  \begin{equation*}
    w_\za=\frac{P(u_1,\ldots,u_n)}{\prod_{i=1}^n u_i^{n_i(\za)}},
  \end{equation*}
  where P is a polynomial such that none of the $u_i$ divides $P$
  $(i=1,\ldots,n)$ and $n_i(\za)$ is the multiplicity of the simple
  module $\za_i$ in the module $\za$.
\end{theorem}
\begin{proof}
  Let $T_C=\{\zb_1,\ldots,\zb_n\}$ be the triangulation of the $(n+3)$
  polygon corresponding to the cluster $C $ and let $D$ be the set of
  diagonals of the polygon; thus $T_C\subset D$. Let
  $T_0=\{-\za_1,\ldots,-\za_n\}$ be the ``snake triangulation''
  \cite[12.2]{cluster2}. ${Q}_{T_0}$ is the alternating quiver of type
  $A_n$ and the diagonals $-\za_i\in T_0$ are the negative simple
  roots. Fomin and Zelevinsky have shown that there is a bijection
  $\za\mapsto w_\za$ between the set of almost positive roots
  $\Phi_{\ge-1}$ and the set of cluster variables $V$. In type $A$,
  they identified $\Phi_{\ge-1}$ with $D $ and proved that for any
  cluster $C$ there is a bijection $\za\mapsto w_\za$ between
  $D\setminus T_C$ and $V\setminus C$. In section \ref{sectionequiv},
  we have shown the bijection $M_\za\mapsto \za$ between $\Ind(\QC) $
  and $D\setminus T_C$. This establishes a bijection $\Ind(\QC)\to
  V\setminus C $. This bijection sends the simple module in
  $\Ind(\QC)$ at the vertex $i\in \QC$ to the variable $w_\zb$ where
  $\zb$ is the unique diagonal in $D\setminus T_C$ that crosses
  $\zb_i$ and does not cross any diagonal in $T_C\setminus\{\zb_i\}$
  
  Let $M_\za$ be an element of $\Ind(\QC)$ with $\za$ the
  corresponding diagonal in $D\setminus T_C$. By Lemma
  \ref{cruciallemma}, we have
  $w_\za=\frac{P(u_1,\ldots,u_n)}{\prod_{i=1}^n u_i^{[\za,\zb_i]}}$.
  We have to show that $[\za,\zb_i]=n_i(\za)$ for all $i=1,\ldots,n$.
  Note that $\za\ne\zb_i$ since $\za\notin T_C$, hence using
  Proposition \ref{prop}, we get $[\za,\zb_i]=(\za\mid\mid\zb_i)$ and
  by \cite[12.2]{cluster2} this is equal to $1$ if the diagonals $\za$
  and $\zb_i$ are crossing, and zero otherwise. Thus
  \begin{equation*}
    [\za,\zb_i]=\left\{
      \begin{array}{ll}
        1&\textup{if $i\in\Supp\za$ in the sense 
          of section \ref{sectionequiv}}\\
        0&\textup{otherwise.}
      \end{array}\right\}
    = n_i(\za).
  \end{equation*}
\end{proof}

\bibliographystyle{alpha}
\bibliography{CaChSc}

\newcommand{\etalchar}[1]{$^{#1}$}
\begin{thebibliography}{MRZ03}

\bibitem[ARS95]{AR}
M.~Auslander, I.~Reiten, and S.~O. Smal{\o}.
\newblock {\em Representation theory of {A}rtin algebras}, volume~36 of {\em
  Cambridge Studies in Advanced Mathematics}.
\newblock Cambridge University Press, Cambridge, 1995.

\bibitem[BFZ]{cluster3}
A.~Berenstein, S.~Fomin, and A.~Zelevinsky.
\newblock {Cluster algebras III: Upper bounds and double Bruhat cells}.

\bibitem[BMR{\etalchar{+}}]{marall}
A.~Buan, R.J. Marsh, M.~Reineke, I.~Reiten, and G.~Todorov.
\newblock Tilting theory and cluster combinatorics.
\newblock In preparation.

\bibitem[FZ02]{cluster1}
S.~Fomin and A.~Zelevinsky.
\newblock Cluster algebras. {I}. {F}oundations.
\newblock {\em J. Amer. Math. Soc.}, 15(2):497--529 (electronic), 2002.

\bibitem[FZ03a]{cluster2}
S.~Fomin and A.~Zelevinsky.
\newblock Cluster algebras. {II}. {F}inite type classification.
\newblock {\em Inventiones Mathematicae}, 154:63--121, 2003.

\bibitem[FZ03b]{Ysystems}
S.~Fomin and A.~Zelevinsky.
\newblock Y-systems and generalized associahedra.
\newblock {\em Annals of Math.}, 158(3), 2003.

\bibitem[Gab80]{gabriel}
P.~Gabriel.
\newblock Auslander-{R}eiten sequences and representation-finite algebras.
\newblock In {\em Representation theory, I (Proc. Workshop, Carleton Univ.,
  Ottawa, Ont., 1979)}, volume 831 of {\em Lecture Notes in Math.}, pages
  1--71. Springer, Berlin, 1980.

\bibitem[MRZ03]{MRZ}
R.~Marsh, M.~Reineke, and A.~Zelevinsky.
\newblock Generalized associahedra via quiver representations.
\newblock {\em Trans. Amer. Math. Soc.}, 355(10):4171--4186 (electronic), 2003.

\end{thebibliography}

\end{document}